\input amstex
\input amsppt.sty
\magnification=\magstep1
\hsize = 6.5 truein
\vsize = 9 truein

\NoBlackBoxes
\TagsAsMath
\catcode`\@=11
\redefine\logo@{}
\catcode`\@=13

 1

    \let    \< = \langle
    \let    \> = \rangle
    \define      \a     {\alpha}
    \redefine   \b      {\beta}
    \redefine   \d      {\delta}
    \redefine   \D      {\Delta}
    \define      \e     {\varepsilon}
    \define      \g     {\gamma}
    \define     \G      {\Gamma}
    \redefine   \l      {\lambda}
    \redefine   \L      {\Lambda}
    \define      \n     {\nabla}
    \redefine   \var    {\varphi}
    \define      \s     {\sigma}
    \redefine   \Sig  {\Sigma}
    \redefine   \t      {\tau}
    \define     \th     {\theta}
    \redefine   \O      {\Omega}
    \redefine   \o      {\omega}
    \define     \z       {\Bbb Z}
    \redefine   \i      {\infty}
    \define      \p     {\partial}

\def\label#1{\par%
        \hangafter 1%
        \hangindent .75 in%
        \noindent%
        \hbox to .75 in{#1\hfill}%
        \ignorespaces%
        }
\document

\baselineskip=16 pt

\centerline{\bf The Kolmogorov-Obukhov Exponent in}
\centerline{{\bf the Inertial Range of
Turbulence: A Reexamination of
Experimental Data}\footnote{Supported in part by the Applied
Mathematical Sciences subprogram of the Office of Energy Research, U.S.
Department of Energy, under contract DE--AC03--76--SF00098, and in part by
the National Science Foundation under grants DMS94--14631 and
DMS89--19074.}}
\bigskip\bigskip
\centerline{\smc G. I. Barenblatt, A. J. Chorin}

\centerline{Department of Mathematics and Lawrence Berkeley Laboratory}
\centerline{University of California}
\centerline{Berkeley, California 94720 USA} \medskip
\centerline{and}
\medskip
\centerline{\smc V. M. Prostokishin}
\centerline{P. P. Shirshov Institute of Oceanology} \centerline{Russian
Academy of Sciences} \centerline{Moscow 117218 Russia}

\vskip .5 truein
\block{\smc Abstract}. In recent  papers
Benzi et al.
presented experimental data and an analysis to the effect
that the well-known ``2/3" Kolmogorov-Obukhov exponent in the inertial
range of local structure in turbulence should be corrected by a small but definitely non-zero
amount. We reexamine the very same data and show that this conclusion is
unjustified.
The data are in fact consistent with incomplete
similarity
in the inertial range, and with an exponent that depends on the Reynolds number and tends to 2/3 in the
limit of vanishing viscosity.
If further data confirm this conclusion, the understanding of local structure
would be profoundly affected.
\endblock

\vskip .5 truein
\specialhead{1. Introduction}\endspecialhead

In 1941 Kolmogorov derived  his famous
scaling relations
for the local structure of turbulence [16]; in particular he deduced that the
second-order
structure function in the inertial range was proportional to the $r^{2/3}$
where $r$ is the separation of the points
(For definitions
and analysis, see below).
Obukhov [18] simultaneously determined the energy spectrum in the
inertial range by similar means.
Soon afterwards, Landau suggested that the Kolmogorov-Obukhov $2/3$
(--5/3)
exponent may
be modified by the presence of intermittency, and various proposals for
Reynolds-number-independent modifications have been made since
then (see e.g. [17]).

On the other hand, the Kolmogorov-Obukhov scaling
argument
has been reexamined
through
the prism of modern scaling theory [4,6,7,9] which produced a
Reynolds-number-dependent
exponent in the structure function, tending to $2/3$ in the limit of
vanishing viscosity.
This argument is supported by the near-equilibrium statistical theory of
turbulence [13,14,15],
as well as by vanishing-viscosity asymptotics [4,5,6,7,14].

The Kolmogorov-Obukhov theory was not derived from first principles such as
as the Navier-Stokes equations, and contains additional assumptions which are
open to debate. In the absence of general analytical solutions of the
Navier-Stokes equations and of adequate computational data, the only way to
check the theory is to subject it to
experimental verification. Several
experimentalists
have claimed to have observed a  correction to the $2/3$ exponent,
and that
it was independent of $Re$; among the influential papers in this
direction are the papers
of Benzi et al. [11,12]. In  [11] (page 389), the authors state that
``the exponents$\dots$are the same
in all experiments" (i.e., they are independent of Reynolds number),
and the
exponent in the second order structure function is ``close but
definitely different from
the value $2/3$ used by Obukhov". Our goal here is to refute
this statement: to the extent that the data  exhibited in [11] can be relied upon,
they militate in favor
of a Reynolds-number
dependent exponent with a $2/3$ vanishing-viscosity limit. The difference
between the
conclusions reached in [11] and the conclusions reached below is apparently
due to
the data not having been carefully enough  examined in [11] for possible dependence
on Re.

In the next section we provide a brief summary of scaling arguments as they
apply
to the structure functions in fully developed turbulence. We then present our
analysis of the data and draw a conclusion.

\vskip .5 truein
\specialhead{2. Scaling in the local structure of turbulence}\endspecialhead

The quantities of interest in the local structure of turbulence are the
moments of the
relative
velocity field, in particular the second order tensor
with components
$$D_{ij}=\<(\Delta_{\bold r})_i(\Delta_{\bold r})_j\>,
\tag 2.1   $$
where $\Delta_{\bold r}={\bold u}({\bold x}+{\bold r})-{\bold
u}({\bold x})$ is a velocity difference between
${\bold x}$ and
${\bold x}+{\bold r}$. In incompressible flow that is in addition locally isotropic, all the
components of this
tensor are
determined if one knows $D_{LL}=\langle [u_L({\bold x}+{\bold
r})-u_L({\bold x})]^2\rangle$ where
$u_L$ is the velocity component along the vector {\bf r}.

To derive an expression for  $D_{LL}$ assume, following
Kolmogorov, that for $r=|\bold r|$ small, it depends
on  $\<\e\>$, the mean rate of energy
dissipation per unit volume, $r$, the distance between the points at
which the velocity is measured, a length scale $\Lambda$,
for example the Taylor macroscale $\Lambda_T$,
and the
kinematic viscosity
$\nu$\,: $$D_{LL}(r)=f(\<\e\>,r,\Lambda_T,\nu), \tag 2.2 $$
where the function $f$ is assumed to be
the same for all developed turbulent flows.
Introduce the Kolmogorov scale $\Lambda_K$, which marks the lower
bound of
the ``inertial" range of scales in which energy dissipation is negligible:
$$\Lambda_K=\frac {\nu^{3/4} } {\<\e\>^{1/4}};\tag 2.3 $$ Clearly,
the  appropriate velocity scale is
$$u=(\<\e\>\Lambda_T)^{1/3},\tag 2.4 $$ and this yields a Reynolds
number
$$Re\ = \ \frac {(\<\e\>\Lambda_T)^{1/3}\Lambda_T}{\nu}\ =\ \frac
{\<\e\>^{1/3}\Lambda_T^{4/3}}{\nu} \ = \
\left(\frac{\Lambda_T}{\Lambda_K}\right)^{4/3}.\tag 2.5 $$

Dimensional analysis yields the scaling law:
$$
D_{LL} =(\langle\e\rangle r)^{\frac 23} \ \Phi
\left(\frac{r}{\Lambda_K} \ ,Re \right) \ , \tag 2.6
$$
where
as before, the function $\Phi$ is an unknown dimensionless function of its
arguments, which have been chosen so that in the inertial range they
are both large.

If one now subjects (2.6) to an assumption of complete
similarity in both its arguments (i.e., one assumes that $\Phi$ tends to a finite non-zero limit when its arguments tend to infinity, see [1]), one obtains
the classical Kolmogorov $2/3$
law [16] $$
D_{LL} =A_0(\<\e\> r)^{\frac 23} \ ,
\tag 2.7
$$
from which the Kolmogorov-Obukhov ``5/3" spectrum [18]
 can be
obtained
via Fourier transform.
If one makes the
assumption of incomplete
similarity in $r/\Lambda_K$ and no similarity in $Re$, (i.e., one assumes that
as $(r/\Lambda_K)\to 0$ the function $\Phi$ has power-type asymptotics with
$Re$-dependent parameters), as we have shown to be appropriate in
certain shear flows [7,8], the result is
$$\frac
{D_{LL}(r)}{(\<\e\>r)^{2/3}} \
= \ C(Re)\left(\frac{r}{\Lambda_K}\right)^{\alpha (Re)},\tag
2.8
$$ where
$C,\alpha$ are functions of $Re$ only.
Expand $C$ and $\alpha$ in powers of $\frac {1}{\ln Re}$,
as is suggested by vanishing viscosity asymptotics, (for another
example of this expansion, see [5,7]),
and keep the two leading terms; this yields

$$ D_{LL}=(\<\e\>r)^{2/3}\left(C_0+\frac{C_1}{\ln
Re}\right)\left(\frac{r}{\Lambda_K}\right)^{\alpha_1/\ln Re},\tag 2.9$$
where $C_0, C_1, \alpha_1$ are constants and the
zero-order term in the exponent has been set
equal to zero  so that $D_{LL}$ has  a finite
limit as
$\nu\to0$ [7].
According to (2.9), the exponent in $D_{LL}$ is  $Re-$dependent
and converges to $2/3$ as $Re\to\infty$.
Note that the prefactor, (the ``Kolmogorov constant"), is also
$Re-$dependent, as has indeed been observed experimentally [19,20].

A further possibility is to subject $D_{LL}$ to an assumption of complete
similarity
in $Re$ and incomplete similarity in $r/\Lambda_K$, opening the door to
$Re$-independent
corrections to the $2/3$ power. This possibility, discussed in [7], is
incompatible
with the existence of a well-behaved vanishing-viscosity limit for the
second-order structure function,
in contradiction with the theory in [5,6,7,14].
This assumption corresponds to the ``extended similarity"
discussed in [11,12].

The conclusion that the classical Kolmogorov-Obukhov value is obtained
in the limit $Re\to\infty$ was reached in [13,14,15] by a statistical mechanics
argument. Furthermore, the usual explanation for
possible departures
of the exponent from the Kolmogorov-Obukhov value is the need to account
for intermittency.
However, it was shown in [13,15] that the Kolmogorov-Obukhov value already
takes intermittency into account; indeed, mean-field theories presented
in [14] give exponents which differ from $2/3$,  and of course depend on the
additional assumptions used to define a system to which
a mean-field theory can be applied. In [7] it was argued that the $2/3$ value
corresponds to ``perfect" intermittency, with the $Re$-dependent
correction being
created by the decrease in intermittency due to viscosity.
This physical picture is mirrored by the fact that we obtained  the $2/3$
exponent
not by an assumption of complete similarity but as the vanishing-viscosity
limit
of a power-law derived from an assumption of incomplete similarity; for a more
detailed explanation in a related problem, see [8].

Kolmogorov [16] proposed similarity relations also for the higher-order
structure functions:
$$D_{LL\dots L}(r)=\langle [u_L({\bold x+\bold r})-
u_L({\bold
x})]^p\rangle,$$ where $LL\dots L$ denotes $L$ repeated $p$ times; the
scaling gives $D_{LL \dots L}=C_p(\langle\e\rangle r)^{p/3}$.
As is well-known, for $p=3$ the Kolmogorov scaling is valid with no
corrections.
We shall not be
concerned here with higher-order structure functions, for which the validity
of the vanishing-viscosity arguments is at present unknown , and whose very
existence in the
limit of vanishing viscosity is doubtful [7].

\vskip .5 truein
\specialhead{3. A reexamination  of the data of Benzi et al.}\endspecialhead

Our starting point is the graph in Figure 3 of the paper [11], which
contains a plot
of the second order structure function $\log D_{LL}$ as a function of the logarithm of the
third moment $D_{LLL}$, whose dependence on $r$  in the inertial interval is well-known to
obey the appropriate Kolmogorov scaling and thus be proportional to $r$. This way of
processing the data provides a longer interval in which the exponent can be seen, and also
produces as an artifact a slope of $2/3$ for separations $r$ in the dissipation
range, as is indeed carefully explained in [11].  The data come from four
series of experiments in a small wind-tunnel: experiments labeled J in which
the turbulence was produced by a jet  at $Re=300,000$ (based on the integral
scale), experiments labeled C6 where the turbulence was produced by a cylinder
at $Re=6000$,
experiments labeled C18, with a cylinder and $Re=18000$, and
experiments in which the
turbulence was produced by a grid and $Re$ was not specified in the paper;
we have found  from referees that the Reynolds number of the grid data was low, with no
precise value. The various
experimental procedures are detailed in [11] and we do not query them in
any way.

The resulting values of $\log D_{LL}$  as a function of $\log D_{LLL}$ were plotted in
Figure 3 of [11]  without regard to $Re$, and a line was fitted to them as
if they
came from a single experiment. That line had slope roughly equal to $0.7$,
and this
is the basis for the claim by Benzi et al. that the exponent is definitely
larger
than $2/3$. However, it is obvious even to the naked eye that the  points
that come
from experiments with differing $Re$ do not lign up well on that single line. We
now show this lack of alignment in detail.

To our regret, Benzi et al. have not responded to our requests for
data in digital
form, so we scanned Figure 3 of [11] with the modern equipment available at the Lawrence Berkeley Laboratory and obtained numerical values in this
way. The set of values is
incomplete because in certain regions of that Figure there are so many
points that it is
impossible to separate them properly; there are however enough scanned points so that the
conclusions below are independent of the remainder.

In Figure 1 we display the  resulting values of $\log D_{LL}$ as a function of
$\log D_{LLL}$ separately for each run, together with lines that correspond to
$y\!=\!(2/3)x+a_1$, $y\!=\!0.7x+a_2$, where
$x, y$ are the coordinates in those graphs and the values of $a_1,a_2$
are the same as in [11].
As one can see,
the slope defined by the experimental points is almost exactly $0.7$ for
the experiment $C6$
($Re=6000$); it goes down as $Re$ increases first to 18,000 and then to 300,000.
The grid data, which we are told belong to a low Reynolds number,  are particularly instructive: They follow the $2/3$ line for
a while, presumably while the separation $r$ is in the
dissipation scale, and then they bend towards the $.7$ line, as
the separation emerges from the dissipation scale but the Reynolds number is
not large enough to produce the asymptotic $2/3$ scaling. The information
at our disposal is not sufficient to estimate {\it a priori} where the bend
should be.
It is quite clear  from these figures that one cannot view all these points as
lying on a single line,
and the data are compatible with incomplete similarity in $r/\Lambda_K$ and an
absence of similarity in $Re$, so that
the exponent in the power law for $D_{LL}$ and therefore
the slopes of the
lines in Figure 1 are slowly decreasing functions of $Re$ when the separation $r$ is in the inertial range.

To make this point another way, we display in Figure 2 the local slopes in
these figures,
defined as $$s^{\text{local}}_i=\frac{y_i-y_1}{x_i-x_1},
\tag3.1$$
where $(x_i,y_i)$ are the coordinates of the i-th point in the graph for
specific experiment,
and $(x_1,y_1)$ is the first (leftmost) point in that graph. The local
slopes that result from using
successive points are too noisy for any conclusion to be drawn.
Figure 2
clearly
shows that the slopes decrease as $Re$ increases for separations $r$
outside the dissipation range. In particular, for the grid data (apparently lowest $Re$)
the slope increases with the separation $r$ as that separation emerges from the
dissipation range.

\vskip .5 truein
\specialhead{4. Conclusion}\endspecialhead

Figures 1 and 2 show, we believe conclusively, that there is no reason to
conclude with [11] that
the Kolmogorov-Obukhov exponent is ``definitely" different from $2/3$ and
independent
of $Re$. The data as we presented them suggest to the contrary that the
exponent
slowly decreases with $Re$ and tends to $2/3$ as $Re\to\infty$. The
experimental uncertainties
detailed in [11], the uncertainty aboout the Reynolds number of the grid data, the small differences between the slopes under
discussion, and
the added uncertainties of the scanning, deter us from making the statement
more emphatic yet.

\newpage
\bigskip\centerline{\bf References}

\roster

\item
G.I. Barenblatt,
{\it Similarity, Self-Similarity and Intermediate Asymptotics}, Consultants
Bureau, NY,
(1979); {\it Scaling, Self-Similarity, and Intermediate Asymptotics},
Cambridge University Press, Cambridge, (1996).

\item
G.I. Barenblatt,
On the scaling laws (incomplete self-similarity with respect to Reynolds
number) in the developed turbulent flow in pipes, {\it C.R. Acad.
Sc. Paris}, series II, {\bf 313}, 307--312 (1991).

\item
G.I. Barenblatt,
Scaling laws for fully developed turbulent shear flows. Part 1: Basic
hypotheses and analysis,
{\it J. Fluid Mech.}, 248, 513--520
(1993).

\item
G. I. Barenblatt and A.J. Chorin, Small viscosity asymptotics
for the inertial range of local structure and for the wall region of
wall-bounded turbulence, {\it Proc. Nat. Acad. Sciences USA} {\bf 93},
6749--6752 (1996).

\item
G.I. Barenblatt and A.J. Chorin, Scaling laws and vanishing viscosity
limits for wall-bounded shear flows and for local structure in developed
turbulence, {\it Comm. Pure Appl. Math.} {\bf 50}, 381--398  (1997).

\item
G.I. Barenblatt and A.J. Chorin, Scaling laws and vanishing viscosity
limits in turbulence theory,
 {\it Proc. Symposia Appl. Math. AMS}, {\bf54}, 1-25,
(1998).

\item
G.I. Barenblatt and A.J. Chorin, New perspectives in turbulence: Scaling laws,
asymptotics and intermittency, in press, SIAM Review (1998).

\item
G.I. Barenblatt, A.J. Chorin, and V.M. Prostokishin, Scaling laws in
fully developed
turbulent pipe flow, {\it Appl. Mech. Rev.}, {\bf50}, 413-429, (1997).

\item
G. I. Barenblatt and N. Goldenfeld,
Does fully developed turbulence exist?
Reynolds number dependence vs.~asymptotic covariance, {\it Phys.~Fluids A}
{\bf7} (12),
3078-3082,  (1995).

\item
G.I. Barenblatt and V.M. Prostokishin, Scaling laws for fully developed
shear flows. Part 2. Processing of experimental data, {\it J. Fluid Mech.}
{\bf 248}, 521--529 (1993).

\item
R. Benzi, C. Ciliberto, C. Baudet and
G. Ruiz Chavarria, On the scaling of three dimensional homogeneous
and isotropic turbulence, {\it Physica D} {\bf 80}, 385-398 (1995).

\item
R. Benzi, S. Ciliberto,
R. Tripiccione, C. Baudet, F. Massaioli, and S. Succi,
Extended self-similarity in turbulent flows, {\it Phys. Rev. E} {\bf 48} (1),
R29--R32, (1993).

\item
A. J. Chorin,
Scaling laws in the vortex lattice model of turbulence,
{\it Commun. Math. Phys.} {\bf114}, 167-176, 1988.

\item
A. J. Chorin, {\it Vorticity and Turbulence}, Springer, 1994.

\item
A. J. Chorin, Turbulence as a near-equilibrium process, {\it Lectures in
Appl. Math.} {\bf 31}, 235--248 (1996).

\item
A. N. Kolmogorov, Local structure of turbulence in incompressible fluid at a
very high Reynolds number, {\it Dokl. Acad. Sc. USSR} {\bf 30},
299-302 (1941).

\item
A. S. Monin and A. M. Yaglom,
{\it Statistical Fluid Mechanics}, Vol.~2, MIT Press, Boston, 1975.

\item
A.M. Obukhov, Spectral energy distribution in turbulent flow, {\it Dokl.
Akad. Nauk USSR},
{\bf1}, 22--24 (1941).

\item
A. Praskovsky and S. Oncley, Measurements of the Kolmogorov constant and
intermittency exponents at very high Reynolds numbers, {\it Phys.~Fluids A} {\bf 6} (9),
2786--2788 (1994).

\item
K. R. Sreenivasan, On the universality of the Kolmogorov constant,
{\it Phys.~A} {\bf 7} (11), 2778--2784 (1995).

\endroster

\newpage
\centerline{\bf Figure Captions}

\bigskip
\noindent Figure 1. Variation of $\log D_{LL}$ with $\log D_{LLL}$ plotted separately
for the several values of $Re$:
\block{ a. Experiment C6 ($Re=6000$); b. Experiment C18 ($Re=18000$); c. Experiment J ($Re=300000$);
d. Grid turbulence.}\endblock

\bigskip\noindent
Figure 2: Local slopes in Figure 1 for the several experiments.
(Diamonds-- C6; squares-- C18; triangles-- J; stars-- grid turbulence).
Note that the slope for grid turbulence start from 2/3 and increases, as one may expect
from the calculation of the dissipation range (see the text).

\enddocument